\documentclass[runningheads]{llncs}
\usepackage[T1]{fontenc}
\usepackage{graphicx}
\usepackage{tikz}
\usepackage{pgfplots}
\usetikzlibrary {arrows.meta} 
\usepackage{color, colortbl}
\usepackage{amsmath}
\usepackage{hyperref}

\newcommand\blfootnote[1]{%
  \begingroup
  \renewcommand\thefootnote{}\footnote{#1}%
  \addtocounter{footnote}{-1}%
  \endgroup
}

\begin{document}
\title{Long-Term Multi-Objective Optimization for Integrated Unit Commitment and Investment Planning for District Heating Networks}
\titlerunning{Multi-Objective Integrated Unit Commitment and Investment Planning}
\author{Stephanie Riedmüller\inst{1}\orcidID{0009-0006-4508-4262} \and
Fabian Rivetta\inst{1}\orcidID{0009-0009-4288-4645} \and
Janina Zittel\inst{1}\orcidID{0000-0002-0731-0314}}
\authorrunning{S. Riedmüller et al.}

\institute{Zuse Institute Berlin, Applied Algorithmic Intelligence Methods Department, Takustraße 7, 14195 Berlin, Germany\\
\email{riedmueller@zib.de, zittel@zib.de} }

\maketitle    

\begin{abstract}

The need to decarbonize the energy system has intensified the focus on district heating networks in urban and suburban areas. Therefore, exploring transformation pathways with reasonable trade-offs between economic viability and environmental goals became necessary. 
We introduce a network-flow-based model class integrating unit commitment and long-term investment planning for multi-energy systems. While the integration of unit commitment and investment planning has been applied to multi-energy systems, a formal introduction and suitability for the application of long-term portfolio planning of an energy provider on an urban scale has yet to be met. Based on mixed integer linear programming, the model bridges the gap between overly detailed industrial modeling tools not designed for computational efficiency at scale and rather abstract academic models. 
The formulation is tested on Berlin's district heating network. 
Hence, the challenge lies in a large number of variables and constraints and the coupling of time steps, for example, through investment decisions.
A case study explores different solutions on the Pareto front defined by optimal trade-offs between minimizing costs and CO$_2$ emissions through a lexicographic optimization approach. 
The resulting solution catalog can provide decision-makers valuable insights into feasible transformation pathways, highlighting distinctions between robust and target-dependent investments.\blfootnote{The content of this report is also available as publication in Operations Research Proceedings 2024. Please always cite as: Riedmüller, S., Rivetta, F., Zittel, J. (2025). Long-Term Multi-objective Optimization for Integrated Unit Commitment and Investment Planning for District Heating Networks. In: Glomb, L. (eds) Operations Research Proceedings 2024. OR 2024. Lecture Notes in Operations Research. Springer, Cham. \href{https://doi.org/10.1007/978-3-031-92575-7_33}{\textcolor{blue}{https://doi.org/10.1007/978-3-031-92575-7\_33}} }

\end{abstract}

\section{Introduction}
While energy providers are urged to decarbonize, current optimization tools for multi-energy systems still fall into two categories: On one hand, detailed industrial modeling tools (e.g,. BoFiT by Volue\footnote{https://www.volue.com/products-and-services/volue-bofit-optimisation-chp-sector-coupling-and-industry}, ResOpt by Kisters\footnote{https://www.kisters.eu/product/forecasting-and-optimisation/}) do not allow for an integrated structural optimization
for large input data due to the inherent complexity. On the other hand, academic models either efficiently solve large-scale systems but lack accuracy (e.g.,  PyPSA \cite{ref_pypsa}, REMix \cite{ref_remix}) or are tailored to small-scale systems (e.g., \cite{ref_integration2}, \cite{ref_ecos2}, \cite{ref_ecos3}). Recent developments in the energy modeling community show an increase in integrating network design and operational decisions 
and apply this approach to micro-grids or short-term scenarios, while a solution for urban areas of large size and long-term scenarios is not available due to efficiency issues.  
Still, no solution effectively bridges the gap between these industrial and academic tools.

We address the research gap of formally defining the problem class of integrated design and operational optimization for energy providers of large urban areas from a mathematical point of view. To that end, we propose a novel model class in Section~\ref{section:model_class} that offers the level of detail and flexibility required for industrial usage while being tractable for long-term planning. We demonstrate its effectiveness by presenting a case study on an instance of Berlin's district heating network, the most complex grid in Western Europe, in Section~\ref{section:case_study}. The study explores a lexicographic bi-objective approach to determine a solution catalog with optimal and relevant trade-offs between cost and environmental targets. We investigate how the conflicting objectives affect the transformation paths and distinguish robust from target-dependent investments.

\section{A Multi-Energy Integrated Operation and Design Optimization Problem} \label{section:model_class}

We introduce a novel network-flow-based model class integrating operation and design optimization for multi-energy systems.
The operation optimization is modeled as a unit commitment problem \cite{ref_uc}, incorporating a time horizon, a set of generating and storage units, a representation of the topology, and demand/price forecast time series as input. Unit commitment involves three decisions: 
\begin{description}
    \item[commitment decisions:] whether a unit is producing/storing energy at a time,
    \item[production decisions:] how much energy a unit is producing/storing at a time,
    \item[network flow decisions:] how much energy is flowing on each edge of the grid at a time. Since the model further includes sources (generation units) and sinks (demand), this part represents a network flow problem.
\end{description}
Design optimization or investment planning comprises: 
\begin{description}
    \item[design/investment decisions:] whether an investment is selected.
\end{description}
Finally, the integrated problem is formulated as a mixed-integer linear program (MILP).

\subsection{MILP} 

Let $R$ be a set of resources, $T$ the set of time steps, $I$ the set of generating units and $K$ the set of storage units.
\subsubsection{Variables:}
For resource $r \in R$ and time $t \in T$, we denote as variables by $x_{t, i^{in}}^r, x_{t, i^{out}}^r \in \mathbf{R}_{\geq 0}$ the flow incoming/outgoing flow at unit $i \in I $, by $z_{i,t} \in \{0,1\}$ the status of unit $i \in I$, by $s_{i,t} \in \{0,1\}$ whether the status of unit $i \in I$ changed at $t$, by $h^r_{t,k}$ the storage level of storage $k \in K$, by $p_t^r, e_t^r \in \mathbf{R}_{\geq 0}$ purchased and sold resource and by $\hat{z}_i \in \{0,1\}$ whether an investment is selected.
\subsubsection{Objectives:}
We minimize for costs $f_1$, which include investment and operational costs containing fuel and power trading, and for CO$_2$ emissions $f_2$.
\subsubsection{Full Program:}
For given demand vector $d$, parameter vectors $a$ and $c$, and for piecewise linearized characteristic curves $\varphi$, the following MILP describes integrated unit commitment and investment planning for multi-energy systems:\\
\resizebox{\linewidth}{!}{
  \begin{minipage}{1.3\linewidth}
  \begin{align}
    &\min && \left(f_1^{inv}(\hat{z}) + \sum_{t \in T} f_1^{op} (z_t, s_t, h_t, p_t, e_t, x_t), \sum_{t \in T} f_2 (p_t, x_t) \right) &&& \nonumber \\
    & && \sum_{i \in I}  x_{t, i^{out}}^{r} + \sum_{k \in K} x_{t, k^{out}}^{r} + p_t^r  = d^r + \sum_{i \in I}  x_{t, i^{in}}^{r} + \sum_{k \in K} x_{t, k^{in}}^{r} + e_t^r &&& \forall r \in R \label{eq:balance}\\
    & && x_{t, i^{out}}^{r_2} = \varphi_{i,t}^{r_1,r_2} \left(x_{t, i^{in}}^{r_1}\right) &&& \forall i \in I, t\in T, r_1, r_2 \in R \label{eq:conversion}\\
    & && s_{i,t} \leq z_{i,t}, s_{i,t} \leq 1 - z_{i,t}, s_{i,t} \geq z_{i,t} - z_{i,t-1} &&& \forall i \in I, t\in T \label{eq:activation}\\
    & && \sum_{\tau \in T^{up}_{i,t}} ( s_{i,t} - z_{i,t}) \leq 0, \sum_{\tau \in T^{down}_{i,t}} ( s_{i,t} + z_{i,t} -1) \leq 0 &&& \forall i \in I, t\in T \label{eq:minupdown}\\
    & && x_{t+1, i^{out}}^{r} - x_{t, i^{out}}^{r} \leq a_i^{up}, x_{t, i^{out}}^{r} - x_{t+1, i^{out}}^{r} \leq a_i^{down} &&&\forall i \in I, t\in T, r \in R \label{eq:ramping}\\
    & && h^r_{t+1,k} = a^{loss}_{t,k}h^r_{t,k} + a^{load}_{t,k}x_{t, k^{in}}^{r} - a^{unload}_{t,k}x_{t, k^{out}}^{r} &&&\forall k \in K, t\in T, r \in R \label{eq:storage}\\
    & && z_{i,t} \leq \hat{z}_i &&&\forall i \in I, t\in T \label{eq:invest}\\
    & && h,x,p,e \leq c_{max}, \qquad h,x,p,e \geq c_{min} \label{eq:capacities} \\
    \nonumber
\end{align}
  \end{minipage}
} \\
\newline
The main goal is to fulfill the given demand while balancing all resources, represented by Constraint~\ref{eq:balance}. Constraint~\ref{eq:conversion} describes the conversion from one resource to another at a generating unit in form of a piecewise linear characteristic curve. Constraints~\ref{eq:activation} - \ref{eq:ramping} incorporate technical constraints such as activation, minimum up and down times, and ramping for generating units. Constraint~\ref{eq:storage} models the time linking of storage levels. Constraint~\ref{eq:invest} links investment decisions to the corresponding operation variables. Capacities are set for all non-binary variables by Constraints~\ref{eq:capacities}.
For an in-depth discussion on the model, we refer to \cite{ref_zibreport}.

\section{Case Study on Berlin's District Heating Network} \label{section:case_study}

\subsection{Instance} 
\begin{figure}
\centering
\resizebox{9.0cm}{!}{
   \begin{tikzpicture}[node distance={30mm},font=\bfseries, 
    thick, main/.style = {draw, circle}, 
    balance/.style = {draw, circle, dotted, purple}, 
    demand/.style = {draw, purple, circle, fill=lightgray}, 
    fuel_node/.style = {draw, circle, fill=olive}, 
    power_node/.style = {draw, circle, fill=teal},
    heat/.style = {-{Stealth[length=3mm]}, purple}, 
    fuel/.style = {-{Stealth[length=3mm]}, olive},
    power/.style = {-{Stealth[length=3mm]}, teal}] 

    \node[demand] (D_M)  {Demand};
    \node[balance] (B1) [right of=D_M] {Balance};
    \node[main] (G_M1) [above right of=B1] {HWE Gas};
    \node[main] (G_M2) [right of=B1] {HWE Gas};
    \node[main] (G_M3) [below right of=B1] {CHP Coal};
    \node[demand] (D_NOR) [below of=D_M] {Demand};
    \node[balance] (B2) [left of=D_NOR] {Balance};
    \node[main] (G_W1) [above left of=B2] {P2H};
    \node[main] (G_W2) [left of=B2] {CHP Coal};
    \node[main] (G_H) [below left of=B2] {CHP Gas};
    \node[balance] (B3) [below of=D_NOR] {Balance};
    \node[main] (G_I) [right of=B3] {HWE Gas};
    \node[demand] (D_SUD) [below of=B3] {Demand};
    \node[balance] (B4) [left of=D_SUD] {Balance};
    \node[main] (G_L1) [left of=B4] {HWE Gas};
    \node[main] (G_L2) [below left of=B4] {CCGT};
    \node[main] (N) [right of=G_I] {$\geq$};
    \node[demand] (D_T) [right of=N] {Demand};
    \node[main] (G_S) [above of=D_T] {HWE Gas};
    \node[demand] (D_K) [right of=D_T] {Demand};
    \node[balance] (B5) [above of=D_K] {Balance};
    \node[main] (G_D1) [above of=B5] {HWE Gas};
    \node[main] (G_D2) [above right of=B5] {CCGT};
    \node[balance] (B6) [below of=D_T] {Balance};
    \node[main] (G_N1) [ left of=B6] {HWE Gas};
    \node[main] (G_N2) [below left of=B6] {CCGT};
    \node[main] (G_P) [below  of=B6] {Heatpump};
    \node[main] (G_T) [below  right of=B6] {HWE Gas};
    \node[main] (G_K) [below of=D_K] {CHP Gas};
    \node[demand] (D_F) [right of=G_K] {Demand};
    \node[fuel_node] (F1) [above left of=G_D1] {Fuel};
    \node[fuel_node] (F2) [below left of=G_W2] {Fuel};
    \node[fuel_node] (F3) [left of=G_N1] {Fuel};
    \node[power_node] (P1) [ left of=G_W2] {Power};
    \node[power_node] (P2) [ below of=D_F] {Power};
    \node[power_node] (P3) [ below of=F3] {Power};
    \node[power_node] (P4) [ right of=G_L2] {Power};
    \node[power_node] (P5) [ below of=G_M3] {Power};
    \node[power_node] (P6) [ below of=G_D2] {Power};

    \draw[heat] (G_M1) to [out=260,in=10,looseness=1] (B1);
    \draw[heat] (G_M2) to [out=180,in=0,looseness=1] (B1);
    \draw[heat] (G_M3) to [out=100,in=350,looseness=1] (B1);
    \draw[heat] (B1) to [out=180,in=0,looseness=1] (D_M);
    \draw[heat] (B1) to [out=240,in=40,looseness=1] (D_NOR);
    \draw[heat] (D_NOR) to [out=60,in=220,looseness=1] (B1);
    \draw[heat] (B2) to [out=0,in=180,looseness=1] (D_NOR);
    \draw[heat] (G_W1) to [out=290,in=160,looseness=1] (B2);
    \draw[heat] (G_W2) to [out=0,in=180,looseness=1] (B2);
    \draw[heat] (G_H) to [out=80,in=200,looseness=1] (B2);
    \draw[heat] (B3) to [out=260,in=100,looseness=1] (D_SUD);
    \draw[heat] (D_SUD) to [out=80,in=280,looseness=1] (B3);
    \draw[heat] (D_NOR) to [out=260,in=100,looseness=1] (B3);
    \draw[heat] (B3) to [out=80,in=280,looseness=1] (D_NOR);
    \draw[heat] (G_I) to [out=180,in=0,looseness=1] (B3);
    \draw[heat] (B4) to [out=0,in=180,looseness=1] (D_SUD);;
    \draw[heat] (G_L1) to [out=0,in=180,looseness=1] (B4);
    \draw[heat] (G_L2) to [out=80,in=200,looseness=1] (B4);
    \draw[heat] (G_S) to [out=270,in=90,looseness=1] (D_T);
    \draw[heat] (G_D1) to [out=270,in=90,looseness=1] (B5);
    \draw[heat] (G_D2) to [out=200,in=70,looseness=1] (B5);
    \draw[heat] (B5) to [out=270,in=90,looseness=1] (D_K);
    \draw[heat] (D_K) to [out=190,in=350,looseness=1] (D_T);
    \draw[heat] (D_T) to [out=10,in=170,looseness=1] (D_K);
    \draw[heat] (D_NOR) to [out=0,in=180,looseness=1] (N);
    \draw[heat] (N) to [out=0,in=180,looseness=1] (D_T);
    \draw[heat] (B6) to [out=90,in=270,looseness=1] (D_T);
    \draw[heat] (G_N1) to [out=330,in=230,looseness=1] (B6);
    \draw[heat] (G_N2) to [out=30,in=250,looseness=1] (B6);
    \draw[heat] (G_P) to [out=90,in=270,looseness=1] (B6);
    \draw[heat] (G_T) to [out=150,in=290,looseness=1] (B6);
    \draw[heat] (G_K) to [out=0,in=180,looseness=1] (D_F);
    \draw[heat] (G_K) to [out=90,in=270,looseness=1] (D_K);
    \draw[fuel] (F1) to [out=180,in=0,looseness=1] (G_M1);
    \draw[fuel] (F1) to [out=200,in=0,looseness=1] (G_M2);
    \draw[fuel] (F1) to [out=220,in=0,looseness=1] (G_M3);
    \draw[fuel] (F1) to [out=240,in=90,looseness=1] (G_S);
    \draw[fuel] (F1) to [out=300,in=150,looseness=1] (G_D1);
    \draw[fuel] (F1) to [out=330,in=90,looseness=1] (G_D2);
    \draw[fuel] (F2) to [out=20,in=270,looseness=1] (G_W2);
    \draw[fuel] (F2) to [out=0,in=180,looseness=1] (G_H);
    \draw[fuel] (F2) to [out=340,in=90,looseness=1] (G_L1);
    \draw[fuel] (F2) to [out=320,in=180,looseness=1] (G_L2);
    \draw[power] (P1) to [out=20,in=180,looseness=1] (G_W1);
    \draw[power] (G_W2) to [out=180,in=0,looseness=1] (P1);
    \draw[power] (G_H) to [out=160,in=340,looseness=1] (P1);
    \draw[fuel] (F3) to  (G_N2);
    \draw[fuel] (F3) to  (G_N1);
    \draw[fuel] (F3) to [out=290,in=200,looseness=1] (G_P);
    \draw[fuel] (F3) to [out=340,in=170,looseness=1] (G_T);
    \draw[fuel] (F3) to [out=50,in=160,looseness=1] (G_K);
    \draw[fuel] (F3) to  (G_I);
    \draw[power] (G_D2) to [out=270,in=90,looseness=1] (P6);
    \draw[power] (G_M3) to [out=270,in=90,looseness=1] (P5);
    \draw[power] (G_K) to  (P2);
    \draw[power] (G_N2) to  (P3);
    \draw[power] (P3) to  (G_P);
    \draw[power] (G_L2) to  (P4);
\end{tikzpicture}}
\caption{Abstract structure of the Berlin instance} \label{fig1}
\end{figure}
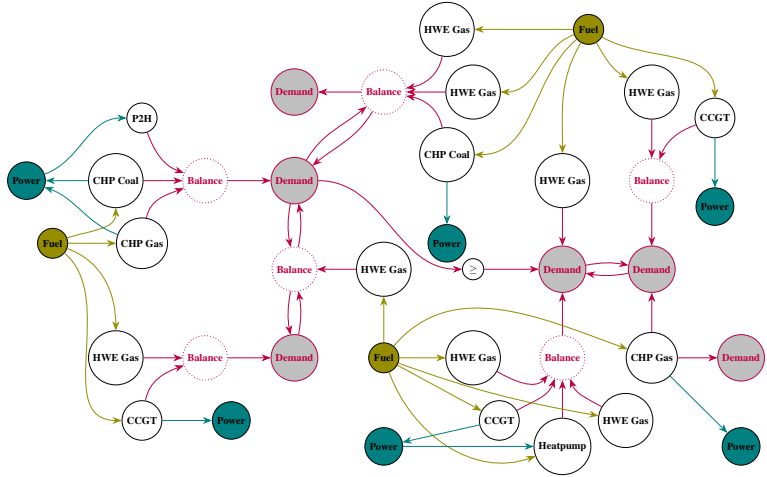

The effectiveness of the model is demonstrated on an instance of Berlin's district heating network. A simplified structure is depicted in Figure~\ref{fig1}.

The initial topology includes heat-generating units with a capacity distribution of 43\% from heating stations and 57\% from combined heat and power plants (CHP), predominantly fueled by gas (75\%) and coal (19\%). In addition, the model incorporates various storage units
and integration with relevant markets. 
Beyond the initial grid topology, the model encompasses 38 
potential investments with a 25-year depreciation period. These include 11 additional heating stations, five storage units, three heat pumps of varying capacities up to 120 MW, four electrical heaters with capacities up to 320 MW, a solar thermal facility, five CHP plants, and nine gas turbines or upgrades to existing gas turbines. 
The resulting MILP contains approximately 3.6M variables and 3.5M constraints.

\subsection{Method} 

The goal is to explore various transformation pathways to achieve decarbonization targets while maintaining reasonable economic costs. Due to the instance size, finding all Pareto optimal solutions (i.e., solutions that cannot be improved in one objective without worsening another) for two objectives (costs, CO$_2$ emissions) is intractable.
Therefore, we explore a reasonable number of Pareto optimal solutions using a lexicographic optimization approach, prioritizing costs before CO$_2$ emissions (hierarchical method \cite{Marler04}). While minimizing CO$_2$ emissions, the optimal cost value is relaxed 
in increments of 5\%, ranging up to 30\% of the cost optimum. 
The 1\% relaxation admits a 2\% MIP-gap, the 5\% relaxation a 1\% MIP-gap, and further relaxations a 0.5\% MIP-gap.
Our computations were executed on a high-performance Cluster, featuring an Intel(R) Xeon(R) Gold 6342 CPU and comparable CPUs with a memory limit of 100 GB. Each individual solution was processed within a time frame of 1 to 6 hours.

\subsection{Computational Results} 

\begin{figure}
\centering
\resizebox{7.9cm}{!}{
    \begin{tikzpicture}
\begin{axis}[
    title={Pareto optimal trade-offs},
    xlabel={normalized costs},
    ylabel={normalized $CO_2$ emissions},
    xmin=0.98, xmax=1.32,
    ymin=0.65, ymax=1.02,
    xtick={1.00,1.05,1.10,1.15,1.20,1.25,1.30},
    ytick={0.5,0.6,0.7,0.8,0.9,1.0},
    xmajorgrids=true,
    ymajorgrids=true,
    grid style=dashed,
    axis x line*=bottom,
    axis y line*=left,
    unit vector ratio*=2 1 1
]
\addplot[color=blue,mark=square,thick,]coordinates {(1.01, 1.0)};
\addplot[color=blue,mark=square,thick,]coordinates {(1.05, 0.90)};
\addplot[color=blue,mark=square,thick,]coordinates {(1.10, 0.82)};
\addplot[color=blue,mark=square,thick,]coordinates {(1.15, 0.78)};
\addplot[color=blue,mark=square,thick,]coordinates {(1.20, 0.73)};
\addplot[color=blue,mark=square,thick,]coordinates {(1.25,0.70)};
\addplot[color=blue,mark=square,thick,]coordinates {(1.30, 0.67)};

\addplot [black, mark = , nodes near coords=\color{blue}{+ electrical heater},every node near coord/.style={anchor=190}] coordinates {(1.25,0.70)};
\addplot [black, mark = , nodes near coords=\color{blue}{- gas turbine},every node near coord/.style={anchor=190}] coordinates {(1.15, 0.78)};
\end{axis}
\end{tikzpicture}}
\caption{Set of computed points on the Pareto front through an iterative relaxation of the costs with corresponding changes in the selected investments.} \label{fig2}
\end{figure}
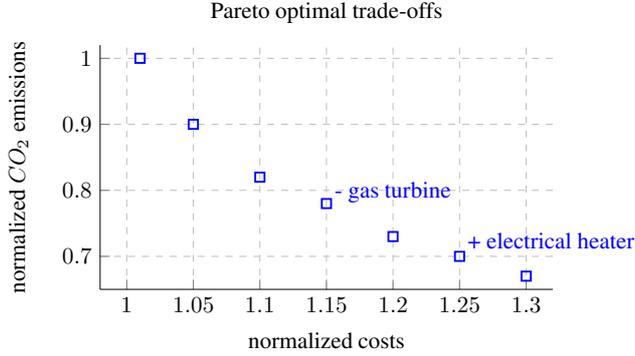

Adapting the cost tolerance gap significantly impacts CO$_2$ emissions. The rate of CO$_2$-savings is particularly rapid for the near-optimal cost solution, gradually diminishing as the tolerance gap widens, as illustrated in Fig.~\ref{fig2}. 

Table~\ref{tab1} presents the resulting investment solution catalogue. Among the potential investments, 27 investments consistently remain unselected, while 10 emerge as fixed choices, i.e. robust investment decisions. 
These include an additional heating station, four new CHP plants, and five gas turbines or their upgrades. 
Three distinct transformation pathways have been identified: The cost-optimal pathway, which includes the robust investments plus an additional gas turbine. Omitting this gas turbine can achieve a CO$_2$-emission reduction of 22-26\%, while introducing an electrical heater results in a CO$_2$-emission reduction of 30-33\%.

\vspace{-2mm}
\begin{table}
\centering
\caption{Selected investments.}\label{tab1}
\resizebox{9.5cm}{!}{
\begin{tabular}{lcccccccl}
\hline
Costs &  101\% & 105\% & 110\% & 115\% & 120\% & 125\% & 130\%\\
\hline
CO$_2$ &  100\% & 90\% & 82\% & 78\% & 73\% & 70\% & 67\%\\
\hline
CHP &  1 & 1 & 1& 1& 1& 1& 1& robust\\
CHP &  1 & 1 & 1& 1& 1& 1& 1& robust\\
Block CHP &  1 & 1 & 1& 1& 1& 1& 1& robust\\
CCGT &  1 & 1 & 1& 1& 1& 1& 1& robust\\
Heating station (biomass) &  1 & 1 & 1& 1& 1& 1& 1& robust\\
Gas turbine upgrade &  1 & 1 & 1& 1& 1& 1& 1& robust\\
Gas turbine &  1 & 1 & 1& 1& 1& 1& 1& robust\\
Gas turbine &  1 & 1 & 1& 1& 1& 1& 1& robust\\
Gas turbine &  1 & 1 & 1& 1& 1& 1& 1& robust\\
Gas turbine &  1 & 1 & 1& 1& 1& 1& 1& robust\\
Gas turbine &  1 & 1 & 1& 0& 0& 0& 0& target-dependent\\
Electrical heater (120 MW) &  0 & 0 & 0& 0& 0& 1& 1& target-dependent\\
\hline
Total &  11 & 11 & 11& 10& 10& 11& 11& \\
\hline
\end{tabular}
}
\end{table}

Notice that the reduction of CO$_2$ emissions is not solely dependent on the choice of investments but is also linked to the operation of a given system setup. As a result, identical investment portfolios can produce different CO$_2$-emission outcomes. For instance, the three scenarios characterized by a cost tolerance gap up to 10\% do not differ in investment decisions but show a CO$_2$-emission reduction of up to 17\%. This highlights the importance of integrating investment planning with unit commitment to make well-informed decisions.

\section{Conclusion}
In conclusion, we formally defined integrated unit commitment and investment planning for multi-energy systems in a mathematical closed form as a MILP. A non-linear formulation of the problem class can be considered by omitting the piecewise linearization of the characteristic curves. In a case study on Berlin's district heating network, we presented a method to efficiently compute bi-objective solutions for the proposed model. 
The findings of this study highlight the interaction of investment and operational decisions. 
Exploring trade-offs within relevant parameter ranges through multi-objective optimization offers a detailed understanding of robust and target-dependent investment decisions and enables stakeholders to make informed decisions for long-term energy portfolio planning.

\begin{credits}
\subsubsection{\ackname} The work for this article has been conducted in the Research Campus MODAL funded by the German Federal Ministry of Education and Research (BMBF) (fund numbers 05M14ZAM, 05M20ZBM). The work has been supported by the German National Science Foundation (DFG) Cluster of Excellence MATH+.

\subsubsection{\discintname}
The authors have no competing interests to declare that are
relevant to the content of this article.
\end{credits}

\bibliographystyle{splncs04}
\bibliography{references}

\begin{thebibliography}{1}
\providecommand{\url}[1]{\texttt{#1}}
\providecommand{\urlprefix}{URL }
\providecommand{\doi}[1]{https://doi.org/#1}

\bibitem{ref_pypsa}
Brown, T., Hörsch, J., Schlachtberger, D.: Pypsa: Python for power system analysis. Journal of Open Research Software  (Jan 2018). \doi{10.5334/jors.188}

\bibitem{ref_integration2}
Castelli, A.F., Pilotti, L., Monchieri, A., Martelli, E.: Optimal design of aggregated energy systems with (n-1) reliability: Milp models and decomposition algorithms. Applied Energy  \textbf{356},  122002 (2024). \doi{10.1016/j.apenergy.2023.122002}

\bibitem{ref_zibreport}
Clarner, J.P., Tawfik, C., Koch, T., Zittel, J.: Network-induced unit commitment - a model class for investment and production portfolio planning for multi-energy systems. ZIB-Report 22-16, ZIB, Takustr. 7, 14195 Berlin (2022)

\bibitem{ref_remix}
Gils, H.C., Scholz, Y., Pregger, T., {Luca de Tena}, D., Heide, D.: Integrated modelling of variable renewable energy-based power supply in europe. Energy  \textbf{123},  173--188 (2017). \doi{10.1016/j.energy.2017.01.115}

\bibitem{ref_ecos3}
Lambert, J., Spliethoff, H.: A two-phase nonlinear optimization method for routing and sizing district heating systems. Energy  \textbf{302},  131843 (2024). \doi{https://doi.org/10.1016/j.energy.2024.131843}

\bibitem{Marler04}
Marler, R., Arora, J.: Survey of multi-objective optimization methods for engineering. Structural and Multidisciplinary Optimization  \textbf{26},  369--395 (04 2004). \doi{10.1007/s00158-003-0368-6}

\bibitem{ref_uc}
Sheble, G.B., Fahd, G.N.: Unit commitment literature synopsis. IEEE Transactions on Power Systems  \textbf{9:1} (Feb 1994). \doi{10.1109/59.317549}

\bibitem{ref_ecos2}
Verheyen, J., Schöbel, N., Thommessen, C., Roes, J., Hoster, H.: Integration of heat pumps and renewable heat suppliers in optimized short-term planning for district heating systems pp. 1419--1430 (2024), the 37th International Conference on Efficiency, Cost, Optimization, Simulation and Environmental Impact of Energy, ECOS 2024

\end{thebibliography}

\end{document}